\documentclass[12pt]{article}
\usepackage{mathrsfs}
\usepackage{amssymb}
\usepackage{amsmath,amsthm}

\newtheorem{Theorem}{Theorem}[section]
\newtheorem{Lemma}[Theorem]{Lemma}
\newtheorem{Corollary}[Theorem]{Corollary}
\newtheorem{Proposition}[Theorem]{Proposition}
\newtheorem{Remark}[Theorem]{Remark}

\theoremstyle{definition}
\newtheorem*{pf}{Proof}

\begin{document}


\title{Number of general Jacobi quartic curves over finite
fields
\thanks{Supported by NSF of China (No. 10990011)}}
\author{Rongquan Feng$^{1}$, Hongfeng Wu$^{2}$
\\ \\
{\small 1~LMAM, School of Mathematical Sciences, Peking
University,}\\{\small Beijing
100871, P.R. China}\\
{\small 2~Academy of Mathematics and Systems Science, Chinese Academy of Sciences,}\\{\small Beijing 100190, P.R. China}\\
{\small\small fengrq@math.pku.edu.cn, whfmath@gmail.com}}

\date{}
\maketitle

\begin{abstract}
In this paper the number of $\overline{\mathbb{F}}_q$-isomorphism
classes of general Jacobi quartic curves, i.e., the number of
general Jacobi quartic curves with distinct $j$-invariants, over the
finite field $\mathbb{F}_q$ is enumerated.

\end{abstract}


\section{Introduction}
Elliptic curve cryptosystems were proposed by Miller (1986) and by
Koblitz (1987) which relies on the difficulty of the elliptic curve
discrete logarithm problem. No sub-exponential algorithms have been
found for solving the discrete logarithm problem based on elliptic
curves is one of main advantage of this system. One basic operation
required to implement the system is the point multiplication, that
is, the computation of $kP$ for an integer $k$ and a point $P$ on
the curve. To obtain faster operations, much effort have been done
in representing the elliptic curves in special forms which provide
faster addition, doubling and tripling in the last decades.

The Jacobi quartic curve is one of the most important curves in
cryptography.  Recent works have shown that arithmetics on the
Jacobi quartic elliptic curves can be performed more efficiently.
The reader is referred to \cite{Bernstein} for the comparison
analysis of computational costs for all kinds of curves. A Jacobi
quartic elliptic curve over a field $K$ is defined by
$y^2=x^4+ax^2+1$, where $a\in K$ with $a^2\neq 4$. Such curves were
first proposed by Chudnovsky and Chudnovsky \cite{Chudnovsky} in
1986. After that, Billet and Joye \cite{Billet}, Duquesne
\cite{Duquesne}, and Hisil, etc. \cite{Hisil} gave more improvements
for the arithmetics on Jacobi quartic curves.

In order to study the elliptic curve cryptosystem, one need first to
answer how many curves there are up to isomorphism, because two
isomorphic elliptic curves are the same in the point of
cryptographic view. In this paper the number of
$\overline{\mathbb{F}}_q$-isomorphism classes of Jacobi quartic
curves, i.e., the number of Jacobi quartic curves with distinct
$j$-invariants, over a finite field is enumerated.

Throughout the paper, $\mathbb{F}_q$ denotes the finite field with
$q$ elements and $\overline{\mathbb{F}}_q$ is the algebraic closure
of $\mathbb{F}_q$.

\section{Background}
A curve means a projective variety of dimension $1$. There are
several ways to define elliptic curves. In this paper, an
irreducible curve is said to be an elliptic curve if it is
birational equivalent to a plane non-singular cubic curve.

It is well-known that every elliptic curve $E$ over a field $K$ can
be written as a Weierstrass equation
$$E:~Y^2+a_1XY+a_3Y=X^3+a_2X^2+a_4X+a_6$$ with coefficients
$a_1,a_2,a_3,a_4,a_6\in K$. The discriminant $\triangle(E)$ and the
$j$-invariant $j(E)$ of $E$ are defined as
$$\triangle(E)=-b_2^2b_8-8b_4^3-27b_6^2+9b_2b_4b_6$$
and
$$j(E)=(b_2^2-24b_4)^3/\triangle(E),$$ where
\begin{equation*}
\begin{array}{rcl}
     b_2 &=& a_1^2+4a_2,\\
     b_4 &=& 2a_4+a_1a_3, \\
     b_6 &=& a_3^2+4a_6,\\
     b_8 &=& a_1^2a_6-a_1a_3a_4+4a_2a_6+a_2a_3^2-a_4^2.
\end{array}
\end{equation*}

Two projective varieties $V_1$ and $V_2$ are isomorphic if there
exist morphisms $\phi: V_1\rightarrow V_2$ and $\varphi:
V_2\rightarrow V_1$, such that $\varphi\circ \phi$ and
$\phi\circ\varphi$ are the identity maps on $V_1$ and $V_2$
respectively. The morphisms $\phi$ and $\varphi$ are called the
isomorphisms from $V_1$ to $V_2$ and from $V_2$ to $V_1$
respectively. Two elliptic curves are said to be isomorphic if they
are isomorphic as projective varieties. Let $$E_1:~
Y^2+a_1XY+a_3Y=X^3+a_2X^2+a_4X+a_6$$ and $$E_2~:
Y^2+a_1^{'}XY+a_3^{'}Y=X^3+a_2^{'}X^2+a_4^{'}X+a_6^{'}$$ be two
elliptic curves defined over $K$. It is known \cite{Silveman} that
$E_1$ and $E_2$ are isomorphic over $\overline{K}$, or $E_1$ is
$\overline{K}$-isomorphic to $E_2$, if and only if $j(E_1)=j(E_2)$,
where $\overline{K}$ is the algebraic closure of $K$. However (see
\cite{Silveman}), $E_1$ and $E_2$ are isomorphic over $\overline{K}$
if and only if there exist $u,r,s,t\in \overline{K}$ and $u\neq 0$
such that the change of variables
$$(X,Y)\rightarrow (u^2X+r,u^3Y+u^2sX+t)$$ maps the equation of $E_1$ to the equation of $E_2$.
Therefore, $E_1$ and $E_2$ are isomorphic over $\overline{K}$ if and
only if there exists $u,r,s,t\in \overline{K}$ and $u\neq 0$ such
that
\begin{equation*}\left\{
\begin{array}{rcl}
     ua_1^{'} &=& a_1+2s,\\
     u^2a_2^{'} &=& a_2-sa_1+3r-s^2, \\
     u^3a_3^{'} &=& a_3+ra_1+2t,\\
     u^4a_4^{'} &=& a_4-sa_3+2ra_2-(t+rs)a_1+3r^2-2st,\\
     u^6a_6^{'} &=& a_6+ra_4+r^2a_2+r^3-ta_3-t^2-rta_1.
\end{array}\right.
\end{equation*}
For the simplified Weierstrass equations where
$a_1=a_3=a_1^{'}=a_3^{'}=0$, then $E_1$ is $\overline{K}$-isomorphic
to $E_2$ if and only if there exist $u,r\in \overline{K}$ and $u\neq
0$ such that
\begin{equation}\label{equiso}
\left\{\begin{array}{rcl}
     u^2a_2^{'} &=& a_2+3r, \\
     u^4a_4^{'} &=& a_4+2ra_2+3r^2,\\
     u^6a_6^{'} &=& a_6+ra_4+r^2a_2+r^3.
\end{array}\right.
\end{equation}
The reader is referred to \cite{Silveman} for more results on the
isomorphism of elliptic curves.

In order to enumerate the number of elliptic curves with distinct
$j$-invariants, one need to study the value distribution of the
$j$-invariant as a function of curve parameters. However it is
effective for only low degree $j$-invariant functions, but very
difficult for counting the value set of the $j$-invariant function
of a Jacobi quartic curve where a polynomial of degree 6 is
involved, especially for the general Jacobi quartic curve where
polynomials with 2 variables are associated. In this paper, this
number is enumerated by studying the
$\overline{\mathbb{F}}_q$-isomorphism classes of those curves.

\section{Enumeration for Jacobi quartics curves}
Let $E_a:~y^2=x^4+ax^2+1$ ($a^2\neq 4$) be a Jacobi quartic curve
defined over a field $K$ of characteristic $>3$. It is clear that
the $j$-invariant of $E_a$ is $\frac{16(a^2+12)^3}{(a^2-4)^2}$.

\begin{Lemma}\label{lem31}
Let $K$ be a field of characteristic $>3$, and let $a\in K$ with
$a^2\neq 4$. Then the curve $$E_a:~y^2=x^4+ax^2+1$$ is birational
equivalent to the elliptic curve
$$W_a:~v^2=u(u-1)\left(u-\frac{2-a}{4}\right)$$ via the change of
variables $\varphi(x,y)=(u,v)$, where
$$u=\frac{x^2-y+1}{2},\quad\mbox{and}\quad v=\frac{x(2x^2-2y+a)}{4}.$$
The inverse change is $\psi(u,v)=(x,y)$, where
$$x=\frac{4v}{4u+a-2},\quad\mbox{and}\quad y=\left(\frac{4v}{4u+a-2}\right)^2-2u+1.$$
\end{Lemma}

\begin{pf}
In order to prove $$v^2=u(u-1)\left(u-\frac{2-a}{4}\right),$$ it is
sufficient to prove $64v^2=4u(4u-4)(4u-(2-a))$. Since
$4u-(2-a)=2x^2-2y+a$ and $64v^2=4x^2(2x^2-2y+a)^2$, it is sufficient
to show that $4x^2(2x^2-2y+a)=4u(4u-4)$. The result then follows
immediately from $$4x^2(2x^2-2y+a)=8x^4-8x^2y+4ax^2,$$
 and
$$\begin{array}{rcl}
4u(4u-4)&=&4(x^2-y+1)(x^2-y-1)\\[.5ex]
&=&4(x^4+y^2-2x^2y-1)\\[.5ex]
&=&4(x^4+x^4+ax^2+1-2x^2y-1)\\[.5ex]
&=&8x^4-8x^2y+4ax^2\end{array}.$$

On the other hand, from $16v^2=4u(u-1)(4u+a-2)$,
$x=\frac{4v}{4u+a-2}$, and
$y=\left(\frac{4v}{4u+a-2}\right)^2-2u+1$, we have $y^2=x^4+ax^2+1$
by a direct computation. Obviously, the maps $\varphi$ and $\psi$
are mutually inverse to each other. \qed
\end{pf}

\begin{Lemma}\label{lem32}
Let $E_a:~y^2=x^4+ax^2+1$ ($a^2\neq 4$) and $E_b:~y^2=x^4+bx^2+1$
($b^2\neq 4$) be two Jacobi quartics curves defined over a field $K$
of characteristic $>3$. Then $j(E_a)=j(E_b)$ if and only if
$\frac{2-b}{4}\in\{\frac{2-a}{4},\frac{4}{2-a},\frac{2+a}{4},\frac{4}{2+a},\frac{2-a}{2+a},
\frac{a+2}{a-2}\}$.
\end{Lemma}

\begin{pf}
The curve $E_a:~y^2=x^4+ax^2+1$ is birational equivalent to the
curve $W_a:y^2=x(x-1)(x-\frac{2-a}{4})$ by Lemma \ref{lem31}.
Therefore $j(E_a)=j(W_a)$. Furthermore, it is well known that for
two Legendre curves $L_{\lambda}: y^2=x(x-1)(x-\lambda)$ and
$L_{\mu}: y^2=x(x-1)(x-\mu)$, they have the same $j$-invariant if
and only if
$\mu\in\{\lambda,\frac{1}{\lambda},1-\lambda,\frac{1}{1-\lambda},\frac{\lambda}{1-\lambda},
\frac{\lambda-1}{\lambda}\}$. Thus the lemma follows. \qed
\end{pf}

\begin{Theorem}\label{jacobi}
Let $N_a$ be the number of $\overline{\mathbb{F}}_q$-isomorphism
classes of Jacobi quartic curves defined over the finite field
$\mathbb{F}_q$. Then we have
\begin{equation*}
N_a=\left\{
\begin{array}{ll}
 \dfrac{q+5}{6},~&\text{if~~} q\equiv 1,7~(\rm{mod}~12),\\[2ex]
 \dfrac{q+1}{6},~&\text{if~~} q\equiv 5,11~(\rm{mod}~12).
\end{array}
\right.
\end{equation*}
\end{Theorem}
\begin{pf} From Lemma \ref{lem32}, we know that
for the elliptic curve $L_{\lambda}: y^2=x(x-1)(x-\lambda)$
($\lambda\neq 0,1$), the map $\lambda \mapsto j(L_{\lambda})$ is
exactly six-to-one unless when $\lambda\in\{-1,2,\frac{1}{2}\}$, the
map is three-to-one, or when $\lambda^2-\lambda+1=0$, the map is
two-to-one. Note that $\lambda^2-\lambda+1=0$ has a root in
$\mathbb{F}_q$ if and only if $\mathbb{F}^{*}_q$ has an element of
order $3$, which is equivalent to $q\equiv 1$ or
$7~(\text{mod}~12)$. Therefore, we have
\begin{equation*}
N_a=\left\{
\begin{array}{ll}
 \dfrac{q-2-3-2}{6}+1+1=\dfrac{q+5}{6},~~&\text{if~} q\equiv
 1,7~(\text{mod}~12),\\[2ex]
 \dfrac{q-2-3}{6}+1=\dfrac{q+1}{6},~~&\text{if~} q\equiv 5,11~(\text{mod}~12).
\end{array}
\right.
\end{equation*}
\qed
\end{pf}

\section{Enumeration for general Jacobi quartics curves}
In this section, consider the general Jacobi quartics curve
$E_{a,b}:~y^2=x^4+ax^2+b$ with $(a^2-4b)b\neq 0$ defined over
$\mathbb{F}_q$ of characteristic $>3$. A Jacobi quartics curve is a
special one of $E_{a,b}$ with $b=1$. The $j$-invariant of $E_{a,b}$
is $j(E_{a,b})=\frac{16(a^2+12b)^3}{b(a^2-4b)^2}$. Note that
$y^2=bx^4+ax^2+1$ can be changed to $y^2=x^4+ax^2+b$ by $x\mapsto
1/x$. So we consider only the form $y^2=x^4+ax^2+b$ for convenience.

The following lemma can be proved by a direct computation similar as
in Lemma \ref{lem31}.

\begin{Lemma}\label{lem41}
Let $K$ be a field of characteristic $>3$, and let $a,b\in K$ with
$(a^2-4b)b\neq 0$. Then the curve $$E_{a,b}:~y^2=x^4+ax^2+b$$ is
birational equivalent to the elliptic curve
$$W_{a,b}:~v^2=u(u^2-2au+a^2-4b)$$ via the change of variables
$$u=2x^2-2y+a,~v=2x(2x^2-2y+a).$$
The inverse change is
$$x=\frac{v}{2u},~y=\left(\frac{v}{2u}\right)^2-\frac{u-a}{2}.$$
\end{Lemma}

For the elliptic curve $E_{a,b}$, we know that
$j(E_{a,b})=\frac{16(a^2+12b)^3}{b(a^2-4b)^2}$. Therefore
$j(E_{a,b})=0$ if and only if $a^2+12b=0$. Moreover, we have the
following proposition.

\begin{Proposition}
Let $E_{a,b}:~y^2=x^4+ax^2+b$ be a general Jacobi quartics curve
defined over the finite field $\mathbb{F}_q$ of characteristic $>3$,
where $(a^2-4b)b\neq 0$. Then $j(E_{a,b})=1728$ if and only if
$a(a^2-36b)=0$, that is, $a=0$ or $a^2=36b$.
\end{Proposition}

\begin{pf}
By Lemma \ref{lem41}, $E_{a,b}$ is birational equivalent to the
curve $W_{a,b}:~y^2=x^3-2ax^2+(a^2-4b)x$, and $W_{a,b}$ is
isomorphic to
$S_{a,b}:~y^2=x^3+(-4a-\frac{a^2}{3})x+(\frac{2a^3}{27}-\frac{8ab}{3})$.
It is clear that the $j$-invariant of $S_{a,b}$ is equal to $1728$
if and only if $\frac{2a^3}{27}-\frac{8ab}{3}=0$, that is
$a(a^2-36b)=0$. Thus $j(E_{a,b})=1728$ if and only if
$a(a^2-36b)=0$.  \qed
\end{pf}

\begin{Corollary} \label{cor}
Let $(a^2-4b)b\ne 0$ and let $N$ be the number of curves
of the form $E_{a,b}$ with $j(E_{a,b})\neq 0,~1728$. Then
\begin{equation*}
N=\left\{
\begin{array}{ll}
 \dfrac{(q-1)(q-7)}{2},~~&\text{if~~} q\equiv 1,7~(\rm{mod}~12),\\
 [2ex]
 \dfrac{(q-1)(q-5)}{2},~~&\text{if~~} q\equiv 5,11~(\rm{mod}~12),
\end{array}
\right.
\end{equation*}
when $b$ is a square and
\begin{equation*}
N=\left\{
\begin{array}{ll}
 \dfrac{(q-1)^2}{2},~~&\text{if~~} q\equiv 1,7~(\rm{mod}~12),\\
 [2ex]
 \dfrac{(q-1)(q-3)}{2},~~&\text{if~~} q\equiv 5,11~(\rm{mod}~12),
\end{array}
\right.
\end{equation*}
when $b$ is not a square.
\end{Corollary}

\begin{pf} Assume first that $b$ is a square in $\mathbb{F}_q$. Then $a^2-4b=0$ has two roots. Hence the number of curves of the
form $E_{a,b}$ over $\mathbb{F}_q$ is
$(q-2)\cdot(q-1)/2=(q-1)(q-2)/2$. If $j(E_{a,b})=0$, then
$a^2+12b=0$ has two roots in $\mathbb{F}_q$ if $q\equiv
1,7~(\text{mod}~12)$, but has no root if $q\equiv
5,11~(\text{mod}~12)$. Therefore the number of curves of the form
$E_{a,b}$ over $\mathbb{F}_q$ with $j(E_{a,b})=0$ is
$2\cdot\frac{q-1}{2}=q-1$ if $q\equiv 1,7~(\text{mod}~12)$, and is 0
if $q\equiv 5,11~(\text{mod}~12)$. If $j(E_{a,b})=1728$, then $a=0$
or $a^2=36b$. Thus the number of curves of the form $E_{a,b}$ with
$j(E_{a,b})=1728$ is
$\frac{q-1}{2}+2\cdot\frac{q-1}{2}=\frac{3(q-1)}{2}$. By
subtraction, we get that
\begin{equation*}
N=\left\{
\begin{array}{ll}
 \frac{(q-1)(q-2)}{2}-(q-1)-\frac{3(q-1)}{2}=\frac{(q-1)(q-7)}{2},~~&\text{if~~} q\equiv
 1,7~(\rm{mod}~12),\\[2ex]
 \frac{(q-1)(q-2)}{2}-0-\frac{3(q-1)}{2}=\frac{(q-1)(q-5)}{2},~~&\text{if~~} q\equiv
 5,11~(\rm{mod}~12).
\end{array}
\right.
\end{equation*}

The number $N$ can be computed similarly when $b$ is not a square.
In this case the number of curves $E_{a,b}$ is
$q\cdot(q-1)/2=q(q-1)/2$ and $a^2+12b=0$ has two roots in
$\mathbb{F}_q$ if $q\equiv 5,11~(\text{mod}~12)$, has no root if
$q\equiv 1,7~(\text{mod}~12)$. \qed
\end{pf}

Now consider curves $E_{a,b}$ over $\mathbb{F}_q$ with
$j(E_{a,b})\ne 0,~1728$. Suppose that two elliptic curves $E_{a,b}$
and $E_{m,n}$ are isomorphic over $\overline{\mathbb{F}}_q$. Then
$j(E_{a,b})=j(E_{m,n})$, and then $j(W_{a,b})=j(W_{m,n})$ by Lemma
\ref{lem41}, which is equivalent to $W_{a,b}$ and $W_{m,n}$ are
isomorphic over $\overline{\mathbb{F}}_q$. Moreover, by
(\ref{equiso}), the last statement holds if and only if there exist
$u,r\in \overline{\mathbb{F}}_q$ with $u\neq 0$ such that
\begin{equation}\label{conreqa}
\left\{\begin{array}{l}
     2m u^2=2a-3r, \\[.5ex]
     (m^2-4n)u^4=3r^2-4ar+(a^2-4b),\\[.5ex]
     r(r^2-2ar+a^2-4b)=0.
\end{array}\right.
\end{equation}

\begin{Proposition}\label{sjacobi}
Let $(a^2-4b)b\ne 0$ and let $b$ be a square element. Then for every
general Jacobi quartic curve $E_{a,b}:~y^2=x^4+ax^2+b$, there is a
Jacobi quartic curve $E_m:~y^2=x^4+mx^2+1$ which is
$\overline{\mathbb{F}}_q$-isomorphic to it.
\end{Proposition}

\begin{pf}
Assume that $b=d^2$ for some $d\in\mathbb{F}^*_q$. Let $m=ad^{-1}$,
$u=d$ and $r=0$. Then $E_{a,b}$ and $E_m$ are isomorphic over
$\overline{\mathbb{F}}_q$ by (\ref{conreqa}).\qed
\end{pf}

\begin{Proposition}\label{prop}
Let $E_{a,b}:~y^2=x^4+ax^2+b$ be a general Jacobi quartics curve
defined over the finite field $\mathbb{F}_q$ of characteristic $>3$,
where $b(a^2-4b)\neq 0$. Assume that $j(E_{a,b})\neq 1728$, then
$E_{a,b}$ and the curve $E_{m,bm^2/a^2}:~y^2=x^4+mx^2+(bm^2/a^2)$
are isomorphic over $\overline{\mathbb{F}}_q$ for any
$m\in\mathbb{F}^{*}_q$.
\end{Proposition}

\begin{pf}
Since $j(E_{a,b})\neq 1728$, we have $a\neq 0$. From
(\ref{conreqa}), for any $m\in\mathbb{F}^{*}_q$, let
$u=\sqrt{\frac{a}{m}}\in\overline{\mathbb{F}}^{*}_q$ and $r=0$, we
get that $E_{a,b}$ and $E_{m,bm^2/a^2}$ are isomorphic over
$\overline{\mathbb{F}}_q$. \qed
\end{pf}

For the elliptic curves $E_{a,b}$ with $j(E_{a,b})\neq 0,\,1728$,
assume that the curve $E_{a,n}$ is
$\overline{\mathbb{F}}_q$-isomorphic to $E_{a,b}$, then there exist
$u,r\in \overline{\mathbb{F}}_q$ with $u\neq 0$ such that
\begin{equation}\label{equ3}
\left\{\begin{array}{l}
     2a u^2=2a-3r, \\[.5ex]
     (a^2-4n)u^4=3r^2-4ar+(a^2-4b),\\[.5ex]
     r(r^2-2ar+a^2-4b)=0.
\end{array}\right.
\end{equation}
just by replacing $m$ to $a$ in (\ref{conreqa}). Thus $r=0$ or
$r^2-2ar+a^2-4b=0$. When $r=0$, we have immediately that $n=b$. In
the following, assume that $r\ne 0$. So we have $r=a+2\sqrt{b}$ or
$r=a-2\sqrt{b}$, Therefore
\begin{equation}\label{nformu1}
     n=\frac{a^2}{4}-\frac{3r^2-4ar+(a^2-4b)}{4u^4}=\frac{a^2(a-2\sqrt{b})^2}{4(a+6\sqrt{b})^2}.
\end{equation}
or
\begin{equation}\label{nformu2}
     n=\frac{a^2(a+2\sqrt{b})^2}{4(a-6\sqrt{b})^2}.
\end{equation}
by substituting $u^2=\frac{2a-3r}{2a}$ and $r=a+2\sqrt{b}$ or
$r=a-2\sqrt{b}$ in the second equation of (\ref{equ3}).

Assume first that $b$ is not a square in $\mathbb{F}_q$. We claim
that neither $\frac{a^2(a-2\sqrt{b})^2}{4(a+6\sqrt{b})^2}$ nor
$\frac{a^2(a+2\sqrt{b})^2}{4(a-6\sqrt{b})^2}$ is an element of
$\mathbb{F}_q$ which is contradictory to $n\in \mathbb{F}_q$. In
fact, if
$$\frac{a^2(a-2\sqrt{b})^2}{4(a+6\sqrt{b})^2}=\frac{a^2}{4(a^2-36b)^2}\cdot
((a-2\sqrt{b})(a-6\sqrt{b}))^2\in \mathbb{F}_q,$$ then
$$((a-2\sqrt{b})(a-6\sqrt{b}))^2=(a^2-8a\sqrt{b}+12b)^2\in\mathbb{F}_q.$$
Therefore, we must have $16a^3+192ab=16a(a^2+12b)=0$. Contradicts to
the assumptions that $j(E_{a,b})\neq 0,~1728$. Similarly, if
$\frac{a^2(a+2\sqrt{b})^2}{4(a-6\sqrt{b})^2}\in \mathbb{F}_q$, one
can get a contradiction again. This proves that if $b$ is not a
square in $\mathbb{F}_q$, and $E_{a,n}$ is
$\overline{\mathbb{F}}_q$-isomorphic to $E_{a,b}$, then we must have
$n=b$. Therefore, by Proposition \ref{prop}, when $b$ is not a
square element, the number of elliptic curves of the form $E_{a,b}$
with $j(E_{a,b})\neq 0,~1728$ in each of its
$\overline{\mathbb{F}}_q$-isomorphism class is $q-1$. Thus when $b$
is not a square, the number of $\overline{\mathbb{F}}_q$-isomorphism
classes of curves of the form $E_{a,b}$ with $j(E_{a,b})\neq
0,~1728$ is
\begin{equation}\label{nonsquare}
\left\{
\begin{array}{ll}
 \dfrac{q-1}{2},~~&\text{if~~} q\equiv 1,7~(\rm{mod}~12),\\[2ex]
 \dfrac{q-3}{2},~~&\text{if~~} q\equiv 5,11~(\rm{mod}~12),
\end{array}
\right.
\end{equation}
by Corollary \ref{cor}.

On the other hand, assume that $b=d^2$ is a square in
$\mathbb{F}_q$. We claim that neither
$\frac{a^2(a-2\sqrt{b})^2}{4(a+6\sqrt{b})^2}$ nor
$\frac{a^2(a+2\sqrt{b})^2}{4(a-6\sqrt{b})^2}$ is equal to $b$. In
fact, if $$\frac{a^2(a-2\sqrt{b})^2}{4(a+6\sqrt{b})^2}=b,$$ then
$$a^2(a-2d)^2=4d^2(a+6d)^2.$$ Thus $$a(a-2d)=2d(a+6d)\;\;\;\;
\mbox{or}\;\; \;\;a(2d-a)=2d(a+6d).$$ So $$(a+2d)(a-6d)=0\;\;\;\;
\mbox{or}\;\;\;\; a^2+12d^2=0,$$ that is $a^2-4b=0$ or $a^2=36b$ or
$a^2+12b=0$, which is contradictory to the assumptions that
$j(E_{a,b})\neq 0,~1728$. Similarly, if
$\frac{a^2(a+2\sqrt{b})^2}{4(a-6\sqrt{b})^2}=b$, then one can get a
contradiction again. Furthermore, we can check easily that
$\frac{a^2(a-2\sqrt{b})^2}{4(a+6\sqrt{b})^2}\neq
\frac{a^2(a+2\sqrt{b})^2}{4(a-6\sqrt{b})^2}$. This proves $n$ has 3
choices in this case. Therefore, when $b$ is a square element, the
number of elliptic curves of the form $E_{a,b}$ with $j(E_{a,b})\neq
0,~1728$ in each of its $\overline{\mathbb{F}}_q$-isomorphism class
is $3(q-1)$ by Proposition \ref{prop}, and then the number of
$\overline{\mathbb{F}}_q$-isomorphism classes of curves of the form
$E_{a,b}$ with $j(E_{a,b})\neq 0,~1728$ is
\begin{equation}\label{square}
\left\{
\begin{array}{ll}
 \dfrac{q-7}{6},~~&\text{if~~} q\equiv 1,7~(\rm{mod}~12),\\[2ex]
 \dfrac{q-5}{6},~~&\text{if~~} q\equiv 5,11~(\rm{mod}~12),
\end{array}
\right.
\end{equation}
when $b$ is a square, by Corollary \ref{cor}.

Adding together the numbers in (\ref{nonsquare}), (\ref{square})
above and 2 which corresponding the two special classes of curves
with $j(E_{a,b})=0$ and $j(E_{a,b})=1728$, respectively, we have the
following theorem.

\begin{Theorem}
Let $N_{a,b}$ be the number of $\overline{\mathbb{F}}_q$-isomorphism
classes of general Jacobi quartic curves $E_{a,b}:~y^2=x^4+ax^2+b$
with $(a^2-4b)b\neq 0$ defined over the finite field $\mathbb{F}_q$.
Then we have
\begin{equation*}
N_{a,b}=\left\{
\begin{array}{ll}
 \dfrac{4q+2}{6},~~&\text{if~~} q\equiv 1,7~(\rm{mod}~12),\\ [2ex]
 \dfrac{4q-2}{6},~~&\text{if~~} q\equiv 5,11~(\rm{mod}~12).
\end{array}
\right.
\end{equation*}
\end{Theorem}

\begin{Remark}
We know from Proposition \ref{sjacobi} that the number $N_a$ of
$\overline{\mathbb{F}}_q$-isomorphism classes of Jacobi quartic
curves $E_a$ is equal to the number of
$\overline{\mathbb{F}}_q$-isomorphism classes of general Jacobi
quartic curves $E_{a,b}$ with square $b$. Thus one have
$$N_a=\frac{q-7}{6}+2=\frac{q+5}{6}$$ when $q\equiv
1,7~(\rm{mod}~12)$ but
$$N_a=\frac{q-5}{6}+1=\frac{q+1}{6}$$ when $q\equiv
5,11~(\rm{mod}~12)$ from the numbers in (\ref{square}) and from
$a^2+12b=0$ has no root, i.e., $j(E_{a,b})\ne 0$, when $b$ is a
square and $q\equiv 5,11~(\rm{mod}~12)$. Therefore we get the result
in Theorem \ref{jacobi} again.
\end{Remark}



\begin{thebibliography}{99}
\bibitem{Bernstein} D.J. Bernstein and T. Lange, Analysis and optimization of
elliptic-curve single-scalar multiplication. Contemp. Math., Vol
461, 1-20, Amer. Math. Soc., 2008.

\bibitem{Billet} O. Billet and M. Joye, The Jacobi model of
an elliptic curve and side-channel analysis, AAECC 2003, LNCS 2643,
34-42, Spriger-Verlag, 2003.

\bibitem{Chudnovsky} D. V. Chudnovsky, and G. V. Chudnovsky, Sequences of numbers
generated by addition in formal groups and new primality and
factorization tests, Advances in Applied Mathematics 7, 385-434,
1986.

\bibitem{Duquesne} S. Duquesne, Improving the arithmetic of elliptic curves
in the Jacobi model, Information Processing Letters 104(3), 101-105,
2007.

\bibitem{Hisil} H. Hisil, G. Carter and E. Dawson, New formulae for
efficient elliptic curve arithmetic, in INDOCRYPT 2007, LNCS 4859,
138-151, Spriger-Verlag, 2007.

\bibitem{Silveman} J.H. Silverman. The Arithmetic of Elliptic Curves, GTM
106, Springer-Verlag, Berlin, 1986.

\end{thebibliography}
\end{document}